\documentclass[12pt, a4paper]{amsart}
\input xypic
\usepackage{amsmath}
\usepackage{amssymb} 
\usepackage{amscd}
\usepackage{latexsym}
\usepackage{amsthm}
\usepackage{amsmath}
\usepackage{amssymb} 
\usepackage{amscd}
\usepackage{latexsym}
\usepackage{amsthm}
\usepackage[all]{xy}

\newcommand{\call}{\mathcal{L}}

\newcommand{\calt}{\mathcal{T}}

\newcommand{\map}{\operatorname{Map}\nolimits}

\def \p{^{\wedge}_p}

\def \nz{\hfill\break}
\def \absatz{\hfill\break\indent}

\def \?{???????????????}

\def \lra{{\, \longrightarrow \,}}

\def \larrow#1{\, \stackrel{#1}\longrightarrow \,} 

\def \Q{\mathbb{Q}}
\def \q{\Q}
\def \Z{\mathbb{Z}}
\def \z{{\Z}}
\def \F{\mathbb{F}}
\def \N{\mathbb{N}}

\def \fp{\F_p}
\def \f2{\F_2}
\def \p{^{\wedge}_p}  %completion
\def \2{^\wedge_2}    %completion
\def \fc2{^\circ_2}  %fiberwise completion
\def \lp{_{(p)}}
\def \l2{{_(2)}}
\def \smash{\wedge}
\def \a7{\minCDarrowwidth{.7cm}} %arrow length for CD
\def \quer{\overline}

\newtheoremstyle{slant}{}{}{\slshape}{}{\bfseries}{.}{.5em}{}%
\newtheoremstyle{special}{}{}{\slshape}{}{\bfseries}{.}{.5em}{\thmnote{#3}}

\newtheorem{theorem}{Theorem}[section]

\newtheorem{proposition}[theorem]{Proposition}
\newtheorem{corollary}[theorem]{Corollary}
\newtheorem{lemma}[theorem]{Lemma}

\newtheorem{remark}[theorem]{Remark}    

\newtheorem{defi}[theorem]{Definition}

\newcommand{\hd}{\operatornamewithlimits{hd}}

\title{Quasi Finite Loop Spaces are Manifolds}

\author{Nitu Kitchloo}
\address{Department of Mathematics, 
404 Krieger Hall, John Hopkins University, Baltimore, MD 21218,
USA}
\email{nitu@math.jhu.edu}
\author{Dietrich Notbohm}
\address{Department of Mathematics \& Computer Science,
University of Leicester, University Road, Leicester, LE1 7RH,
U.K.} 
\email{dn8@mcs.le.ac.uk}


\begin{document}
\begin{abstract}It is an old conjecture, that finite $H$-spaces are 
homotopy equivalent to manifolds. Here we prove that this conjecture
is true for loop spaces. Actually, we show that every quasi finite 
loop space
is equivalent to 
a stably parallelizable manifold. The proof is conceptual and relies on 
the theory of p-compact groups.
 On the way we also give a complete 
classification of all simple 2-compact groups of rank 2.

\end{abstract}

\maketitle

\date{\today}

\def \absatz{\hfill\break\indent}

% macros

%
%
%
%%
%
%

%
%
%%
%
%
%
%
%
%
\hfill\number\day /\number\month /\number\year
%\nz\nz
%
%
%
%
%
%
%
%
%
\begin{section}{Introduction}

It is an old question in the theory of $H$-spaces,
whether finite $H$-spaces are equivalent to differentiable
manifolds. The first major result in this direction is due to Browder 
 who showed in a series of papers 
 that every simply connected finite $H$-space is 
homotopy equivalent to a closed topological manifold and, if the dimension is
not congruent to 2 mod $4$, then this manifold can be taken to be smooth and 
stably parallelizable. 

The first examples of finite $H$-spaces which are
not compact Lie groups, were constructed using Zabrodsky's method of 
mixing homotopy types \cite{hiltonroitberg} \cite{zabrodsky}.
Pedersen analyzed Zabrodsky's method in detail and showed
that, in particular, $H$-spaces in the genus of a compact Lie group
are  homotopy 
equivalent to stably parallelizable
smooth manifolds \cite{pedersena} \cite{pedersenb}. 

In \cite{capellweinberger}, Capell and Weinberger got further results
of this type for finite $H$-spaces. They put some extra conditions
on the fundamental group; e.g that the fundamental group is a 
finite p-group (p odd) 
or infinite with at most 2-torsion. Under these assumptions
they were able to show that such finite $H$-spaces are equivalent to 
topological manifolds.

In this paper we will concentrate on finite loop spaces in general,
and show that all finite loop spaces are homotopy equivalent to
stably parallelizable smooth manifolds.

A loop space  is a triple $(L, BL,e)$ where $L$ and $BL$ are topological
spaces, with $BL$ pointed, and where $e : \Omega BL \lra L$ is a 
homotopy equivalence. By abuse of notation we denote this loop space also
by $L$. Then $L$ is an $H$-space with classifying space $BL$.
Properties of loop spaces are inherited from $L$; e.g 
$L$ is called finite if the space $L$ is homotopy equivalent
to a finite $CW$-complex, 
quasi finite 
if $H^*(L;\z)$ vanishes in large degrees and is a finitely 
generated abelian group in each degree,
and simply connected if the space $L$ is so. 
Since all components of an $H$-space are homotopy equivalent, we can restrict 
ourselves to connected loop spaces. 

\begin{theorem}\label{maintheorem}
For any connected quasi finite loop space $L$, there exists
a stably parallelizable, smooth, 
 finite dimensional closed manifold $M$ 
such that $L$ and $M$ are homotopy equivalent.
\end{theorem}

This theorem also says that every quasi finite loop space 
is actually finite. This was already proved in \cite{nofinob} by methods
similar to those also used in this paper.

The proof of the theorem is based 
on ideas and techniques developed by 
Pedersen in the above mentioned work \cite{pedersena} \cite{pedersenb}.
He constructed a special 1-torus for the $H$-spaces $X$ he considered.
A special 1-torus is a fibration of the form $S^1 \lra X \lra Y$ with 
special extra
properties (see the next section). In particular, the map
$S^1 \lra X$ factors through the inclusion $S^1 \subset S^3$.
Using the fact that
$Y$ is a quasi finite, stably reducible, nilpotent Poincar\'e complex 
he could prove
that $X$ is homotopy equivalent to a stably parallelizable manifold.

The advantage of working with finite loop spaces comes 
from the fact that
after p-adic completion we get p-compact groups. Hence we can treat
$L\p$ almost like a compact Lie group. This will allow us to
construct particular subgroups of $L\p$,
which are used to construct special 1-tori for $L$.

The theory of p-compact groups will also 
provide
enough information to show that $L/S^1$ is a quasi finite, stably reducible,
nilpotent Poincar\'e complex, which is the other main ingredient to
make Pedersen proof applicable. Our proof relies on work 
of Bauer \cite{bauer} who used ideas of Klein \cite{Klein} to construct 
an analogue of the one point compactification of
the adjoint representation of a compact Lie group for p-compact groups
(see Section 3).

Unfortunately, special 1-tori do not exist for all quasi finite loop spaces;
e.g. they do not exist for products of $SO(3)$'s. 
The rational cohomology $H^*(L;\q)$ of a quasi finite loop space is 
an exterior algebra generated by odd dimensional classes. 
We say that $L$ is small if $H^*(L;\q)$ is generated by 
classes  of degree less than or equal to 3. Otherwise, we call $L$ 
large. We have to treat small and large quasi finite loop spaces 
differently. For large loop spaces we follow the ideas of Pedersen and
construct special 1-tori. For small quasi finite loop spaces
we have a complete classification.

\begin{theorem} \label{smallcase} 
Let $L$ be a small, connected, $\z$-finite loop space.
Then there exists a compact Lie group $G$ isomorphic to a product of
$S^3$'s, a central elementary abelian subgroup $E\subset G$, and a torus $T$
such that $G/E\times T$ and $L$ are homotopy equivalent.
\end{theorem}

\begin{remark}
{\rm 
Since any finite loop space $L$ (actually any finite $H$-space)
 is equivalent to a product $L'\times T$ 
where $T$ is a torus and $L'$ is a finite loop space (respectively, a finite 
$H$-space) with finite 
fundamental group we may assume for both theorems 
that $\pi_1(L)$ is finite.
And this we will do in all what follows. We call a quasi finite loop space 
semi simple if it is connected and if $\pi_1(L)$ is finite.
Hence we only need to prove the above theorems 
for semi simple quasi-finite loop spaces. 
}
\end{remark}

The paper is organized as follows.
In the next section we prove Theorem 1.2.
Section 3 is devoted to a 
 discussion of stable reducibility and  completions.
In Section 4 we discuss the notion of special 1-tori and reduce the 
proof of our main theorem for large finite loop spaces 
to the existence of local special 1-tori. 
In  Section 5 we classify simple 2-compact groups of rank 2.
This is needed to
construct completed special 1-tori for 
p-compact groups. This  is worked out in Section 6.
In Section 7 we use an arithmetic square argument to 
construct special 1-tori for 
finite loop spaces. 

We had developed our own argument for proving stable reducibility of $L/S^1$
but unfortunately we couldn't apply it in our situation (see Section 8).
We explain our argument in the last section, since we think 
that it is interesting in it's own right and since it shows the power
of Klein's homotopy theoretic version of the adjoint representation of 
a compact Lie group \cite{Klein}. The techniques are similar to those used in 
\cite{bauer}, but much simpler and very much motivated by
a geometric argument for the stable reducibility of homogeneous spaces.

Several proofs of this paper rely on the theory of p-compact groups.
For a general reference we refer the reader to the survey articles
\cite{mollersurvey} and \cite{nosurvey} and the references given there.

We thank Erik Pedersen for bringing this question  to the 
attention of the first author, for valuable discussions and 
his continuous interest
in this work. We also thank the CRM in Barcelona for it's support
when part of this work was done.

\end{section}
\bigskip
\begin{section}{The case of small finite loop spaces}

In this section we assume that $L$ is a  small, semi simple, connected, 
quasi finite loop space.
We notice that p-adic completion makes $L$ into a p-compact group.
That is $(L\p,BL\p,e\p)$ is a connected p-compact group.

Theorem 1.2 is part of the following statement.

\begin{theorem}
Let $L$ be a connected quasi finite loop space.  Then the following 
conditions are 
equivalent:
\nz
(i) The exterior algebra $H^*(L;\q)$ is generated by 
3-dimensional classes; i.e. $L$ is small.
\nz
(ii)
The polynomial algebra $H^*(BL;\Q)$ is generated by  
4-dimensional classes.
\nz
(iii) There exists a compact connected Lie group $G$ isomorphic to a product 
of  $S^3$'s and  a central elementary abelian subgroup $E\subset G$ 
such that $L\simeq G/E$ (as spaces). 
\end{theorem}

\begin{proof}
In the case of a $p$-compact group, Theorems of this type are proved in
\cite{dwabelian} (Theorems 0.5A and 0.5B). We will be using those results.

The equivalence of (i) and (ii)  follows easily from an Eilenberg-Moore 
spectral sequence argument. 

Let us assume that (ii) holds, and let $L\p$ denote the associated 
p-compact group.
If $p=2$ we can apply Theorem 0.5B of \cite{dwabelian}.
In our case this says that there exists a connected compact Lie group $G$
isomorphic to a product of $S^3$'s and a central elementary abelian 
subgroup $E\subset G$ such that $BL\2 \simeq B(G/E)\2$.
Obviously, the number of $S^3$'s is determined by the number of 4-dimensional
generators of $H^*(BL\2;\z\2)\otimes \q$ or of $H^*(BL;\q)$.

Now we consider the case of an odd prime p.
The universal cover of $\widetilde L\p$ of $L\p$ is 
a simply connected p-compact group \cite{mncenter} and splits into a 
product of simply connected simple p-compact groups
\cite{nosplitting} or \cite{dwsplittings}, i.e. 
$\widetilde L\cong \prod _i X_i$.
Since $H^*(B\widetilde L\p;\z\p)\otimes \q \cong H^*(BL\p;\z\p)\otimes \q$,
the number of factors is determined by the number of generators of
 $H^*(BL;\q)$
and, for each factor $X_i$,
the algebra $H^*(BX_i;\z\p)\otimes\q$ is a 
monogenic polynomial algebra generated by a 4-dimensional class. 
Hence $W_{X_i}\cong \z/2$ and $W_{L\p}\cong \prod_i W_{X_i}$ is an
elementary abelian 2-group. Now
Theorem 0.5A of \cite{dwabelian} tells us that 
$B\widetilde L\p\simeq BH\p$
where $H$ is a product of $S^3$. Since $H\p$ is center free,
$BH\p\simeq B\widetilde L\p= BL\p$. The number of factors again
equal the number of generators of $H^*(BL;\q)$. 

By construction $G\cong H$,
and since $B(G/E)\p\simeq BG\p$ for odd primes, we have 
$BL\p\simeq B(G/E)\p$
for all primes. Since $BL$ and $BG/E$ are both rationally 
products of Eilenberg-MacLane spaces and 
since they have the same rational cohomology, we also know that the 
rationalizations 
$BL_0$ and $B(G/E)_0$ are equivalent. That is to say that $BL$ is 
in the adic genus of $BG/E$ and hence $L$ in the adic genus of $G/E$. 
By \cite{nofinob} (Theorem 6.1), the adic genus of $G/E$ is rigid,
that is every space in the adic genus of $G/E$ is actually 
equivalent to $G/E$.
This implies condition (iii).
 
On the other hand if $L\simeq G/E$, then the exterior algebra 
$H^*(L;\q)$ is generated by 3-dimensional classes. This finishes the proof.
\end{proof}

\end{section}

\bigskip

\begin{section}{Stable reducibility and completions}

For a quasi-finite space $K$ we define the  homological dimension 
$\hd_\z (K)$ of $K$ by the degree of the largest non vanishing 
integral homology group.
 For a Poincar\'e complex this equal the formal 
dimension 
of $K$ which is given by the degree of the 
fundamental class. 
We call a quasi-finite Poincar\'e complex  $K$ of homological dimension 
$\hd_\z(X)=n$
stably reducible, if, for some $r\in\N$ there exists
a map $S^{n+r} \lra S^r\smash K=\Sigma^r K$ such that
$H_{n+r}(S^{n+r};\z) \lra H_{n+r}(\Sigma^r K;\z)$ is an isomorphism  
or, equivalently, if  
the Hurewicz map $h:\pi_{n+r}(\Sigma^r K) \lra H_{n+r}(K;\z)$
is an epimorphism. 
These conditions are also  equivalent to the 
fact that the Spivak normal bundle is stably trivial and to the fact that 
the top cell splits off stably.

Using the techniques of completions we will break the question of 
stable reducibility down to local ones. First we have to recall some notions.

Let $R$ be a commutative ring with unit. 
A space $K$ is called $R$-finite if $H_*(K;R)$ vanishes 
in large degrees and is a finitely generated $R$-module in each degree.
In particular, quasi finiteness is nothing but $\z$-finiteness. 
For such spaces the R-homological dimension of $K$, denoted by 
$\hd_R(K)$, is given by th degree of the largest non vanishing homology group 
(with coefficients in the ring $R$).  
A $R$-finite space $K$ 
with $\hd_{R}(K)=n$ is a 
Poincar\'e complex if $H^*(K;R)$ and $H_*(K;R)$ satisfy the usual
Poincar\'e duality properties with respect to a fundamental class
$[K]_R \in H_n(K;R)$.

We call a $R$-finite Poincar\'e complex $K$ with 
$\hd_R(K)=n$ R-stably 
reducible if, for some $r\in\N$, there exists
a map $S^{n+r} \lra \Sigma^r K$ such that the induced map 
\nz
$H_{n+r}(S^{n+r};R) \lra H_{n+r}(\Sigma^r K;R)$ is an isomorphism 
or, equivalently, if the Hurewicz map
\nz
$\pi_{n+r}(\Sigma^rK)\otimes R \lra H_{n+r}(\Sigma^rK;R)$ is an epimorphism.

For a space $K$ we denote the p-adic completion by $K\p$. If $K$ is nilpotent
or $p$-good 
completion induces an isomorphism in mod-p homology and cohomology.
Hence, for such spaces $K$ and $K\p$ have the same mod-p properties.

\begin{lemma}
Let $K$ be a $\z$-finite, nilpotent, Poincar\'e complex of 
homological dimension $n$.
Then the following are equivalent:
\nz
(i) $K$ is stably reducible.
\nz
(ii)
For all primes $p$,   $K$ is mod-p stably reducible. 
\nz
(iii)
For all primes $p$, the completion $K\p$ is mod-p stably reducible.
\end{lemma}

\begin{proof}
Since $K$ is a Poincar\'e complex of dimension $n$, we know that
$H_{n+r}(\Sigma^rK;\fp)\cong H_{n+r}(\Sigma^rK;\z)\otimes \fp\cong \fp$.
We consider the exact sequence
$$
\pi_{n+r}(\Sigma^r K) \lra H_{n+r}(\Sigma^rK;\z)\cong \z \lra Q \lra 1
$$ 
where $Q$ is the 
cokernel of the 
Hurewicz map. We choose $r$ big enough so that we are in the stable range
and $Q$ does not depend on $r$.
Then  $Q=0$ if and only if $Q\otimes \fp=0$ for all primes. That is 
$K$ is stably reducible if and only if $K$ is mod-p stably reducible 
for all p.

The equivalence between (ii) and (iii) follows from the fact that 
$\pi_*(\Sigma^r K)\otimes \fp \cong \pi_*(\Sigma^r K\p)\otimes \fp$.
\end{proof}

We record another easy lemma.

\begin{lemma}\label{localsplit}
Let $K$ be a $\z\lp$-finite, nilpotent, Poincar\'e complex of 
$\z\lp$-homological dimension $n$. Then the following are equivalent:
\nz
(i) $K$ is $\z\lp$ stably reducible.
\nz
(ii)
$K$ is mod-p stably reducible. 
\nz
(iii)
$K\p$ is mod-p stably reducible.
\end{lemma}

\begin{proof}
The equivalence of (ii) and (iii) is exactly as in the above lemma. 
Hence one only needs to establish the equivalence of (i) and (ii). 
This follows rather trivially since the map 
$\pi_{n+r}(\Sigma^r K)\otimes \z\lp \lra H_{n+r}(\Sigma^rK;\z\lp)$ is 
an epimorphism if and only if the map 
$\pi_{n+r}(\Sigma^r K)\otimes \fp \lra H_{n+r}(\Sigma^rK;\fp)$ is an 
epimorphism. 
\end{proof}

In \cite{bauer} (Theorem 1.3), Bauer showed that for a p-compact group $X$
there exists a p-completed sphere spectrum $S_X$ 
with a stable $X$-action
whose dimension equals 
the $\fp$-homological dimension $\hd_{\fp}(X)$.
Here we have to work in the category of simplicial spaces and have to replace
the loop space $X$ by the associated simplicial group.
Moreover, if $X$ is the completion of a compact Lie group, we can take 
for $S_X$ the one point 
compactification of the Lie algebra of $G$  with the adjoint action.
In particular, if $G$ is abelian, $S_G$ has the trivial action.
He also showed that, for a p-compact subgroup $Y \subset X$ there
exists a map of spectra $S_X \lra X\smash_Y S_Y$ which induces an isomorphism
in $H_n(-;\fp)$. Here, we take the mod-p cohomology of spectra.
This result enables us to prove the following proposition.

\begin{proposition}\label{X/Tred}
Let $X$ be a p-compact group and let $T \lra X$ be a p-compact subtorus.
Then, the $\fp$-finite homogeneous space $X/T$ is mod-p stably reducible.
\end{proposition}

\begin{proof}
Since $T$ is the completion of an abelian compact Lie group, 
the sphere spectra $S_T$ carries the trivial $T$ action. 
Therefore, we get a map $S_X \lra X_+\smash_T S_T \simeq X/T_+ \smash S_T$
which induces an isomorphism in mod-p homology in the right degree.
In particular, this tells us that the top cell of $X/T$ splits off 
stably and that $X/T$ is mod-p stably reducible.
\end{proof}

\end{section}
\bigskip
\begin{section}{Special 1-tori and the proof of Theorem 1.1}

In \cite{pedersenb} Pedersen introduced the concept of special 1-tori 
for spaces,
which  is his main concept to get control of the surgery obstructions
(see \cite {pedersenb} Proposition 2.1). We will recall his notion, 
Actually, we only need the p-local version.
A fibration $F \lra E \lra B$ is 
called orientable if 
$\pi_1(B)$ acts trivially on the set $[F,F]$ of homotopy classes  of 
self equivalences
of $F$.

\begin{defi}
A nilpotent space $K$ admits a p-local special 1-torus if, up to homotopy, 
there exists a diagram
of orientable fibrations 
$$
\diagram
S^1\lp \rto \dto & S^3\lp \rto \dto & S^2\lp \dto \\
S^1\lp \rto \dto & K \rto \dto & B \dto \\
* \rto & A \rto & A
\enddiagram
$$
 such that 
\nz 
(i) $A$ is $\z\lp$-finite. 
\nz
(ii) $B$ is $\z\lp$-finite and $\z\lp$-stably reducible.
\nz
(iii) Localized at 0, the diagram is homotopy equivalent to
$$
\diagram
S^1_0 \rto \dto & S^3_0 \rto \dto & S^2_0 \dto \\
S^1_0 \rto \dto & A_0\times S^3_0 \rto \dto & A_0\times S^2_0 \dto\\
* \rto & A_0 \rto & A_0
\enddiagram
$$
where all vertical fibrations are trivial.
\end{defi}

Using the notion of special 1-tori Pedersen could prove the following
result (see \cite{pedersenb} Theorem 1.4).

\begin{theorem} (Pedersen)\label{pedersenmain}
Let $X$ be a $\z$-finite $H$-space. If for every prime p the 
localization $X\lp$ admits a special 1-torus, then $X$ is homotopy 
equivalent to 
a smooth stably parallelizable manifold.
\end{theorem}

Obviously, there also exists a notion of a global special 1-torus. 
For the proof of the above theorem, Pedersen first showed that, 
under the above assumption,  $X$ has a global special 1-torus. Then he 
used this extra structure to prove that the finiteness obstruction for $X$
vanishes. 
Since $H$-spaces are stably reducible 
\cite{browderspanier}, their Spivak normal bundle is stably trivial.
The existence of a special 1-torus then also implies the vanishing of the 
surgery obstruction for the existence of a homotopy equivalence to 
a stably parallelizable manifold
(see \cite{pedersenb} 2.1 and 4.2).

To prove Theorem \ref{maintheorem},
 it is  therefore only left to show that any large 
$\z$-finite 
loop space  admits p-locally a special 1-torus.
And this is a consequence of the following proposition.

\begin{proposition}\label{S3subgroup}
Let $L$ be a large connected $\z$-finite loop space. Then there exists a 
loop space $N$ and 
a fibration $A \lra BS^3\lp \larrow{f} BN\lp$ such that 
$A$ is simple, $\z\lp$-finite and such that 
$N$ and $L$ are homotopy equivalent spaces. 
Moreover, localized at 0, 
there exists a left inverse $s : BM_0 \lra BS^3_0$ of $f$, 
i.e. $sf_0=id_{BS^3_0}$.
\end{proposition}

The proof of this proposition will be given in Section 7.

\begin{corollary} \label{plocaltorus}
For a large $\z$-finite loop space $L$ and a prime $p$, the localization 
$L\lp$ admits a p-local special 1-torus.
\end{corollary}

\begin{proof}
Let $S^1 \subset S^3$ be the maximal torus of $S^3$. 
Since the loop space $N$ of the last proposition is equivalent to $L$ 
we only have to prove the claim for $N$ or equivalently, we may assume that
there exist a fibration $BS^3\lp \lra BL\lp$ with the desired properties.

 Passing to classifying spaces and localizations, and taking homotopy fibers
we get a commutative diagram of fibration sequences
$$
\diagram
 & & & & A \dto \\
S^1\lp \rto \dto^= & S^3\lp \rto \dto & S^2\lp \rto \dto 
   & BS^1\lp \rto^i \dto^= & BS^3\lp \dto^f \\
S^1\lp \rto \dto^= & L\lp \rto \dto & B \rto \dto 
   & BS^1\lp \rto^g  & BL\lp  \\ 
* \rto & A \rto^= & A
\enddiagram
$$
Here $B$ is the homotopy fiber of the composition  
$BS^1\lp \lra BS^3\lp \lra BL\lp$. As the homotopy fiber of maps between 
simply connected spaces, $A$ and $B$ are simple. 

The three left columns of the above diagram will 
establish a p-local special 1-torus for $L\lp$.
All rows of this $3\times 3$-diagram are  given by
principal fibrations and therefore 
orientable. The same holds for the two left columns.
For the right column we have a pull back diagram
$$
\diagram
S^2\lp \rto \dto^= & B \rto \dto & A \dto \\
S^2\lp \rto & BS^1\lp  \rto  & BS^3\lp  
\enddiagram
$$
The bottom row is an orientable fibration. Hence, 
this also holds for the top row.
This shows that the above $3\times 3$-diagram 
consists of orientable fibrations.

Since $A$ is $\z\lp$-finite, a Serre spectral sequence argument shows that 
the same holds for $B$. 

Localized at 0, there exists a left inverse $s : BL_0 \lra BS^3_0$.
Since $sg_0=sf_0i_0=i_0$, this left inverse establishes 
rationally compatible left inverses 
for all vertical arrows between the second and third row of the above
large  diagram.
In particular this shows that, localized at 0, the  vertical fibrations 
of the 
$3\times 3$-diagram are trivial and that this diagram satisfies 
the third condition of special 1-tori.

To complete the proof it remains to show that $B$ is $\z\lp$-stably 
reducible. We pass to completions. Then $L\p$ becomes a p-compact group.
We get a fibration
$B\p \lra BS^1{}\p \lra BL\p$. Since $B$ was $\z\lp$-finite and simple,
$B$ and $B\p$ have isomorphic mod-p homology. This shows that 
$B\p$ is $\fp$-finite, that $S^1{}\p \lra L\p$ is a monomorphism
of p-compact groups and that $B$ is equivalent to the homogeneous space 
$L\p/S^1{}\p$. By Proposition \ref{X/Tred}, $B$ is $\fp$-stably reducible
and by Lemma \ref{localsplit}, $B$ is $\z\lp$-stably reducible.
This completes the proof and shows that $L\lp$ admits a p-local special 
1-torus.  
\end{proof}

{\it Proof of Theorem \ref{maintheorem}}:
We already discussed th ecase of small 
quasi finite loop spaces.
 Let $L$ be a large $\z$-finite loop space. By Corollary \ref{plocaltorus}
every localization $L\lp$ admits a p-local special 1-torus. 
By Theorem \ref{pedersenmain} this implies that $L$ is 
homotopy equivalent to a smooth stably parallelizable manifold.
\qed

\begin{remark}
{\rm
Unfortunately, there exists no global version of Proposition \ref{S3subgroup};
i.e. a large quasi finite loop space $N$ might not contain a $S^3$ or a $S^1$
as a subgroup. Hence, in general there exists no sequence of fibrations 
of the form $$S^1 \lra L \lra L/S^1 \lra BS^1 \lra BN,$$
where $L$ and $N$ are homotopy equivalent.
If such a sequence were to exist, the argument of Section 8 would give a proof of the stable reducibility of $L/S^1$.
}
\end{remark}

\end{section}
\bigskip
\begin{section}{2-compact groups of rank 2}

In this section we will classify all simple 2-compact groups of rank 2.
For the convenience of the reader and to fix notation we recall some 
material about p-compact groups.

A p-compact group $X$ is a loop space $X=(X,BX,e)$ such that $BX$ is 
$p$-complete 
and pointed and such
that $X$ is $\fp$-finite. 
Every p-compact group $X$ has a maximal torus $T_X$,
 a maximal 
torus normalizer $N_X$, and a Weyl group $W_X$
acting on $T_X$. These loop spaces fit into
a diagram
$$
\diagram
BT_X \rto \drto & BN_X \rto \dto & BW_X \\
& BX
\enddiagram
$$
Here, $BT_X\simeq K((\z\p)^n,2)$ is homotopy equivalent to 
an Eilenberg-MacLane space of degree 2.
The top row is a fibration and determines the action of $W_X$ on $T_X$,
actually on $L_X:=\pi_1(T_X)\cong (\z\p)^n$. We call $L_X$ the associated 
$W_X$-lattice and $n$ the rank of $X$.
This action can also be described by a representation $W_X \lra Gl(L_X)$.
If $X$ is connected, 
this representation is faithful and 
makes the finite group $W_X$ into a pseudo reflection group. 
And if in addition $p=2$,  then $W_X$ is 
a 2-adic reflection group. 
We call $X$ simple if $X$ is connected and if the associated representation
$W_X \lra Gl(L_X\otimes \q)$ is irreducible.
For details and further notions 
we refer the reader to the survey articles \cite{mollersurvey} and 
\cite{nosurvey} and the references mentioned there.

The following theorem might be  known to the experts. But since we 
couldn't find a 
reference for it, we will also include a proof.

\begin{theorem}\label{classification}
Any simple 2-compact group $X$ of rank 2 is isomorphic 
to the 2-adic completion of 
$SU(3)$, $Spin(5)=Sp(2)$,
$SO(5)$ or $G_2$.
\end{theorem}

The rest of this section is devoted to the proof of this statement.
For compact connected Lie groups we will abuse 
notation and denote by $G$ the associated 
2-compact group obtained by 2-completion.

Let $U$ be a finite dimensional vector space over $\q\2$ with an action 
of a finite group $W$ defined by a homomorphism $W\lra Gl(U)$.
A $W$-lattice $L$ of $U$ is a $\z\2$-lattice $L\subset U$ of maximal rank
fixed under the action of $W$; i.e. $L$ is a $\z\2[W]$-module
and $L\otimes \q \cong U$. 
We say that two $W$-lattices $L$ and  $L'$ of $U$ are isomorphic 
if $L\cong L'$ as $\z\2[W]$-modules.
A $W_1$-lattice $L_1$ and a $W_2$-lattice $L_2$ are called weakly isomorphic
if there exists an isomorphism $W_1\cong W_2$ such that $L_1$ and $L_2$
are isomorphic as $W_1$ lattices.

We say that  two p-compact groups $X$ and $Y$ have the same 
Weyl group data if the representations $W_X\lra Gl(L_X)$ and 
$W_Y \lra Gl(L_Y)$ are
weakly isomorphic.  Renaming the elements of $W_Y$ we always can 
identify $W_Y$ with $W_X$ and assume that the two lattices are actually 
isomorphic.

From the Clark-Ewing list \cite{clarkewing} we get a complete list
of all irreducible  reflection groups of rank 2
defined over $\q\2$. These are given by the dihedral groups 
$D_6$, $D_8$ and $D_{12}$ with their standard representation as  
reflection groups. In fact, these are the only dihedral groups 
which can be represented as reflection groups over $\q\2$.
 The first is the rational Weyl group representation of $SU(3)$, the second 
of $Spin(5)$ or $SO(5)$ and the last of the exceptional Lie group $G_2$.
The classification of Clark and Ewing only works up to weak equivalence. 

The Lie groups $Spin(5)$ and $Sp(2)$ are isomorphic. Hence, the
Weyl groups $W_{Spin(5)}$ and $W_{Sp(2)}$ are also isomorphic 
and the
associated lattices $L_{Spin(5)}$ and $L_{Sp(2)}$ weakly isomorphic.
In the following, we will always use the one of these two which seems to be
more appropriate. 

The universal cover $\widetilde X$ of a p-compact group $X$ is again a 
p-compact group and, if $\pi_1(X)$ is finite, $X$ and 
$\widetilde X$ have the same rational Weyl group data and 
$\widetilde X\cong \widetilde X/Z$
where $Z \subset \widetilde X$ is a central subgroup \cite{mncenter}. 
Simple p-compact groups have finite fundamental groups \cite{mncenter}.
Therefore, Theorem \ref{classification} is a consequence of the 
following classification result for simply connected simple 2-compact groups.

\begin{theorem} \label{scclassification}
Let $G=SU(3),\ Sp(2)$ or $G_2$. A simply connected 
2-compact group $X$ has the  same 
rational Weyl group data as $G$ if and only if $X$ and $G$ are 
isomorphic as 2-compact groups.
\end{theorem}

For the proof of this theorem
we first have to classify all 2-adic lattices of the representation 
$W_G \lra GL(L_G\otimes \Q)$.

\begin{lemma}\label{lattices} 
Let $U:=(\q\2)^2$ and $W\lra GL(U)$ be 
a reflection group.
\nz
(i)
If $W=D_6,\ D_{12}$ then, up to isomorphism,
there exists exactly one $W$-lattice of the representation
$W\lra Gl(U)$.
\nz
(ii)
If $W=D_8$ each $W$-lattice of $U$ is isomorphic either to
$L_{SO(5)}$ or to $L_{Spin(5)}$. And both lattices are weakly isomorphic.
\end{lemma}

\begin{proof}
For $D_6$ and $D_{12}$ this follows from \cite{andersen} (Proposition 4.3
and Theorem 6.2).

Now let $W=D_8$. In this case we have two non isomorphic lattices 
$L_{SO(5)}$ and $L_{Spin(5)}$.
Let $L$ be another $W$-lattice of $U$.
For a large $r$, the lattice $2^rL$, the
submodule of all elements divisible by $2^r$, is a submodule of 
$L_{SO(5)}\cong \z\2\oplus \z\2$.
We choose $r$ minimal with this property, i.e. $2^rL\subset L_{SO(5)}$
but $2^{r}L\not\subset 2L_{SO(5)}$. Since 
$L\otimes \q\cong L_{SO(5)}\otimes \q$,
we get a short  exact sequence of $\z\2[W]$-modules 
$$
0\lra 2^rL \lra L_{SO(5)} \larrow{\rho} Q \lra 0.
$$
The minimality of $r$ implies that $Q$ is a finite cyclic group; 
i.e. $Q\cong \z/2^s$ generated by either $\rho((1,0))$ or 
$\rho((0,1))$. The dihedral group $D_8$ is generated by the three elements 
$\sigma_1,\sigma_2, \tau$, where $\sigma_i$ multiplies the i-th coordinate by 
$-1$ and $\tau$ exchanges the two coordinates.
Since the automorphism group of $Q$ is abelian, the action of $W$ on $Q$ factors through the abelianization of $W$, $ab(W)$. It follows that the element $\sigma_1\sigma_2 = \sigma_1\tau\sigma_1\tau$ acts trivially on $Q$. Hence the elements 
$(1,0), (0,1)\in M$
are mapped onto elements of order 2 in $Q$. Thus, either $Q=0$ or $Q=\z/2$.
In the first case, we have $L\cong M=L_{SO(5)}$. In the second case,
$D_8$ acts trivially on $Q$ with $\rho((1,0))=\rho((0,1))\neq 0$ in $\z/2$
and consequently $L\cong L_{Spin(5)}$.
This proves the first part of (ii). 

The second part follows from the facts that 
$L_{Sp(2)}$ and $L_{Spin(5)}$ are weakly isomorphic 
and that $L_{Sp(2)}$ and $L_{SO(5)}$ are isomorphic.   
\end{proof}

{\it Proof of Theorem \ref{scclassification}:}
If $X$ and $G$ have the same rational Weyl group 
data, the above lemma shows that they also have the same 2-adic 
Weyl group data. We can assume that $W:=W_G=W_X$ and that 
$L:=L_G=L_X$. We also can identify the maximal tori $T:=T_G\cong T_X$.

For $G=SU(3)$ or $G_2$, this implies  
$X\cong G$.
For $SU(3)$ this follows
from \cite{mntori}
and for $G_2$ from \cite{viruel}.

For $Spin(5)=Sp(2)$ uniqueness result are only known in terms of the maximal 
torus normalizer \cite{noo(n)} \cite{vavpeticviruel}. We have to show
that $N_X \cong N_{Sp(2)}$ as loop spaces; i.e. $BN_X\simeq BN_{Sp(2)}$.

Since $X$ and $Sp(2)$ have the same rational Weyl group data, they have 
isomorphic rational cohomology. Hence,
$H^*(X;\z\2)\otimes \q$ is an exterior algebra with generators in degree
3 and 7.
If $H^*(X;\z\2)$ has 2-torsion, then $X$ and $G_2$
have isomorphic mod-2 cohomology 
\cite{hubbuck}. The Bockstein spectral sequence then shows that 
$X$ does not have the right rational cohomology.
Therefore, $X$ has no 2-torsion,
$H^*(X;\z\2)$ is an exterior algebra with generators in degree
3 and 7 and $H^*(BX;\f2)\cong \f2[x_4,x_8]$
is a polynomial algebra generated by a class of degree 4 and one of degree 8.
Since $H^*(BX;\f2)$ is a finitely generated module over $H^*(BT;\f2)$, 
the composition 
$$
H^*(BX;\f2)\cong H^*(BX;\z\2)\otimes \f2 \lra H^*(BT;\z\2)^W \otimes \f2
\cong H^*(BSp(2);\f2)
\lra H^*(BT;\f2)
$$ 
is a monomorphism. The isomorphism 
$H^*(BT;\z\2)^W \otimes \f2
\cong H^*(BSp(2);\f2)$
follows from the fact that $X$ and $Sp(2)$
have the same 2-adic
Weyl group data (Lemma \ref{lattices}).
Since the first and third term are both polynomial algebras
of the same type, 
$$
H^*(BX;\f2)\lra H^*(BT;\z\2)^W \otimes \f2\cong H^*(BSp(2);\f2)
$$
is an isomorphism. 

Let $t\subset T$ denote the elements of 
of order 2 and $H:=S^3\times S^3 \subset Sp(2)$ the obvious subgroup.
We have a chain of inclusion $t\subset T \subset H \subset Sp(2)$
and $H=C_{Sp(2)}(t)$. The action of $D_8$ on $t$ factors through 
the $\z/2$-action on $t$ given by switching the coordinates.

Now we use Lannes' T-functor theory (e.g. see \cite{schwartz}).
We get a map $f : Bt \lra BX$ which looks in mod-2 cohomology like the map
$Bt \lra BSp(2)$. This map is $\z/2$-equivariant up to homotopy.
The mod-2 cohomology of 
the classifying space $BC_X(t):=map(Bt,BX)_f$ of the centralizer $C_X(t)$
can be calculated with the help of Lannes' $T$-functor 
and $H^*(BC_X(t);\f2) \cong H^*(BC_{Sp(2)}(t);\f2) \cong H^*(BH;\f2)$.
Moreover, the Weyl group of $C_X(t)$ is given by
the elements of $D_8$ acting trivially on $t$. Hence 
$W_{C_X(t)}\cong \z/2 \times \z/2$.
By \cite{dwabelian} (Theorem 0.5B), this implies that $BC_X(t)\simeq BH$.
We will identify $C_X(t)$ with $H$.
The $\z/2$-action on $t$ induces a $\z/2$-action on $H$.
Since $Bt \lra BX$ was $\z/2$-equivariant up to homotopy,
the inclusion $BC_X(t) \lra BX$ extends to a map
$BY:=BH_{h\z/2} \lra BX$. In this case, the homotopy orbit space 
$BY$ happens to be a 2-compact group and has the same Weyl group as 
$X$. That is $N_Y=N_X$. Moreover, the space $BY$ fits into a fibration
$$
BH \lra BY \lra B\z/2\ ,
$$
which is classified by obstructions in 
$H^{*}(B\z/2;\pi_*(BSHE(BH)))$. Here, $SHE(BH)$ is the space of 
self equivalences of $BH$ homotopic to the identity.
Since $SHE(BH)\simeq (B\z/2)^2$ \cite{dwcenter} and  
since $\z/2$ acts on $\pi_2(B^2(\z/2)^2)\cong (\z/2)^2$ 
by switching the coordinates, all obstruction groups vanish 
and the above fibration
splits. This shows that $BY\simeq B(H\rtimes \z/2):=BH'$ and that
$BN_X = BN_Y \simeq BN_{H'}= BN_{Sp(2)}$.  That is $X$ and $Sp(2)=Spin(5)$ 
have
isomorphic maximal torus normalizer and shows that $X\cong Sp(2)$.
\qed

\begin{remark}
{\rm
The only simply connected 2-compact group of rank 1 is $S^3$.
Hence, we get the following complete list (up to isomorphism)
of connected 2-compact groups of 
rank 2, namely $S^1\times S^1$, $S^1\times S^3$, $U(2)$, $S^1\times SO(3)$,
$S^3 \times S^3$, $S^3\times SO(3)$, $SO(4)$, $SU(3)$, $Sp(2)$,
$SO(5)$ and $G_2$.
}
\end{remark}

The following corollary is needed for later purpose. 
For a p-compact group $X$ we denote by $\quer X$ 
the associated center free quotient.

\begin{corollary}\label{subgroup}
For any simple connected 2-compact group $X$ of rank 2, there
exists a homomorphism $S^3 \lra X$ such that the composition 
$S^3 \lra X \lra \quer X$ is a monomorphism.
\end{corollary}

\begin{proof}
Because of Theorem \ref{classification}
we only have to check this for the compact connected Lie groups 
$SU(3)$, $Sp(2)$, $SO(5)$ and $G_2$.
There exists a chain of monomorphisms
$S^3=SU(2) \subset SU(3) \subset G_2$. Both groups, $SU(3)$ and $G_2$, 
are 2-adically center free. This proves the claim in these two cases.
Let $S^3 \subset Sp(2)$ denote the inclusion into the first coordinate.
Since the intersection of $S^3$ and the center of $Sp(2)$ is trivial,
the composition $S^3 \subset Sp(2) \lra SO(5)$ is also a monomorphism.
This proves the claim in the other cases.
\end{proof}
\end{section}
\bigskip

\begin{section}{Particular subgroups of $p$-compact groups.}

In this section we will construct particular subgroups of large p-compact 
subgroups.
A p-compact group $X$ is called large, 
if the exterior algebra $H^*(X;\z\p)\otimes \Q$
has a generator of degree
$\geq 5$. We want to prove the following proposition.

\begin{proposition}\label{pcgsubgroup}
Let $X$ be a large semi simple connected p-compact group.
Let $r:=dim_{\q\p}\, H^4(BX;\z\p)\otimes \Q$ be the dimension of the
$\q\p$-vector space $H^4(BX;\z\p)\otimes \Q$.
Then there exists a compact Lie Group $G$ and a map 
$$
f: BG\p \lra  BX
$$
such that the following hold:
\nz
(i) $G\cong S^3\times H$ with $H$ semi simple and it's universal cover 
$\widetilde H$ isomorphic to $(S^3)^{r-1}$. If $p$ is odd, 
we can choose $G=(S^3)^{r}$.
\nz
(ii)
The induced map 
$H^4(BX;\z\p)\otimes \q \lra H^4(BG\p;\z\p)\otimes \q$ is an isomorphism.
\nz
(iii)
The homotopy fiber $X/G\p$ of $f$ is simple and $\fp$-finite.
\end{proposition}

\begin{proof}
Comparing the statement with Proposition 3.1 of \cite{nofinob}
there is an extra assumption on the generators of $H^*(X;\z\p)\otimes \q$ 
and the corresponding additional output is that $G$ contains a factor $S^3$.
Actually, for odd primes, the statements of both propositions are the same.
Therefore we only have to prove the statement for $p=2$.
Again, for a compact connected Lie group we denote by $G$ the associated 
2-compact group.

Let $X$ be a semi simple 2-compact group, i.e. $\pi_1(X)$ is
a finite 2-group. 
If $H^*(X;\z\2)\otimes \q$ is not generated by classes of degree 3, i.e.
the polynomial algebra $H^*(BX;\z\2)\otimes \q$ is not 
generated by classes of degree 4, then the 
Weyl group $W_{X}$ is non abelian \cite{dwabelian} (Theorem 0.5B), 
but a honest reflection 
group, since $W_{X}$ is defined over $\q\2$. That is, 
$W_{X}$ is generated 
by elements of order 2 fixing a hyperplane of codimension 1.
The universal cover $\widetilde X$ of $X$ splits into a direct product
$\widetilde X \cong \prod X_i$ 
of simple, simply connected pieces \cite{dwsplittings}. 
Since $X$ and $\widetilde X$ have isomorphic 
Weyl groups, we can assume that $X_1$ has  a non abelian Weyl group $W_1$. 
For this piece we will construct a monomorphism
$BG_1:=BS^3 \lra BX_1$, 
such that the composition 
$BS^3 \lra BX_1 \lra B\quer X_1$ is also a monomorphism (see below). 
Here, 
$\quer X_1$ denotes the center free quotient of $X_1$.
Moreover, the map $BG_1 \lra BX_1$ 
will induce an isomorphism in $H^4(-;\z\p)\otimes \q$.  

Having done this we can proceed similarly 
as in \cite{nofinob}.
For all other pieces there exists monomorphisms $BG_i\lra BX_i$
inducing an isomorphism on $H^4(-;\z\p)\otimes \q$ such that $G_i$ is 
isomorphic 
to $S^3$ or to $SO(3)$ (see \cite{nofinob}). This produces a homomorphism 
$\prod G_i \lra \prod X_i \cong \widetilde X \lra X$
of p-compact groups.
The kernel $K$ of this homomorphism, which might be nontrivial, is a 
central subgroup of $G_1 \times \prod_{i>1} G_i$ \cite{nokernel}.
Since the center free quotient $\quer X$ is isomorphic to 
$\prod_i \quer X_i$ we have a homomorphism $X \lra \quer X_1$.
By construction the composition $S^3 \lra X_1 \lra \quer X_1$ is 
a monomorphism.
We get a commutative diagram
$$
\diagram
K \rto \dto^= & S^3\times  \prod_{i>1} G_i \rto \dto & X \dto \\
K \rto & S^3 \rto & \quer X_1
\enddiagram
$$
where the right arrow in the bottom row is a monomorphism.
Since $\quer X_1$ is center free  
the composition $K \lra S^3 \times \prod_{i>1} G_i \lra S^3$
is trivial. Therefore, $K$ is a subgroup of $\prod_{i>1} G_i$
and the map $S^3 \times \prod_{i>1} G_i \lra X$ factors 
through a monomorphism
$G:=S^3 \times ((\prod _{i>1} G_i)/K) \lra X$
 with all the desired properties.

  It remains to show that, for a simple, simply connected 
$2$-compact group $X$ with non abelian Weyl group, there exists
a monomorphism $S^3 \lra X$ inducing a isomorphisms in
$H^4(-;\z\p)\otimes \q$ such that $S^3 \lra X \lra \quer X$ is also
a monomorphism.
Let $W' \subset W_X$ be a subgroup of the Weyl 
group of $X$ generated by two non commuting reflections of $W_X$.
 Let $T\subset T_X^{W'}\subset T_X$ denote the connected 
component of the fixed-point set of the $W_X$-action on $T_X$, which  
has codimension 2. 
The centralizer $C:=C_X(T)$ is a connected 2-compact group, whose 
Weyl group $W_C$ contains $W'$ \cite{mncenter}. There exists a 
finite covering of 
$C$ which splits into a product $Y\times T$ where $Y$ is a simply connected 
2-compact group of rank 2 with Weyl group isomorphic to $W_C$. 
The action of $W'$ on the maximal torus $T_Y$ of $Y$ gives rise to an
irreducible representation over $\q\2$. Otherwise, $W'$ would split into 
a product and the two chosen reflections would commute.
Hence, the 2-compact group $Y$
is simple and of rank 2. 
By Corollary \ref{subgroup} there exists a monomorphism
$S^3 \lra Y$ such that $S^3 \lra Y \lra \quer Y$ is a monomorphism.
Putting all these homomorphisms  and groups into a diagram
we get
$$
\diagram S^3 \rto \dto^= & Y \rto \dto & X \dto \\
S^3 \rto & Y/K \rto & \quer X
& \ .
\enddiagram
$$
Here, $K$ denotes the kernel of $Y \lra \quer X$. In particular, $K$ is 
a central subgroup of $Y$.
Since $S^3 \lra Y \lra Y/K \lra\quer Y$ is  monomorphism, this also holds
for the composition of the first two arrows. Moreover,  since 
$Y/K \lra \quer X$ is a monomorphism, the same holds for the composition 
$S^3 \lra Y/K \lra \quer X$. This proves  the above claim and finishes 
the proof 
of the proposition.
\end{proof}

\end{section}

\bigskip

\begin{section}{Proof of Proposition \ref{S3subgroup}}

In this section, we want to prove Proposition \ref{S3subgroup}.
The proof is based on  an arithmetic square argument. First we need a statement about the existence 
of a particular 
sub loop space, a global version of Proposition \ref{pcgsubgroup}.

\begin{proposition}\label{globalsubgroup}
Let $L$ be a large semi simple $\z$-finite loop space.
Then there exists a semi simple compact Lie group $G$,
loop spaces $M$ and $N$ and a fibration
$$
A \lra BM \lra BN
$$
such that the following hold:
\nz
(i) $A$ is simple and $\z$-finite.
\nz
(ii)
$G\cong S^3\times H$ and the universal cover of $H$ is isomorphic 
to a product of $S^3$'s.
\nz
(iii)
The spaces $G$ and $M$ as well as $L$ and $N$ are homotopy equivalent.
\nz
(iv) $H^4(BN;\q) \lra H^4(BM;\q)$ is an isomorphism.
\nz
(v) There exists a commutative diagram
$$
\diagram
BM\p \rto \dto & BN\p \dto \\
BG\p \rto & BL\p 
\enddiagram
$$
where the vertical maps are equivalences.
The same holds for the rationalizations of the classifying spaces.
\end{proposition}

\begin{proof}
This statement is a refinement of Proposition 1.4 of \cite{nofinob}.
The proof of that statement is an arithmetic square
argument which uses it's p-completed version, 
namely Proposition 3.1
of \cite{nofinob}, as input. The proof  
carries over word for word. 
We only have to replace that proposition by a
p-completed version of the above claim, namely by
Proposition \ref{pcgsubgroup}.
In particular, the bottom rowin the diagram of
(v) is the map constructed in Proposition \ref{pcgsubgroup}.
Claim (ii), which is not part of Proposition 1.4 of \cite{nofinob},
is a consequence of the same formula in Proposition \ref{pcgsubgroup}.
\end{proof}

\begin{remark}
{\rm
The above proposition establishes an oriented fibration
$G\lra L \lra L/G$. And the existence of such an oriented fibration
is already sufficient to show that the finiteness obstruction  vanishes and 
that
every quasi finite loop space is actually finite (see \cite{nofinob}).
The existence of a special tori is needed for the vanishing of 
the appropriate surgery obstruction.
}
\end{remark}

For the proof of Proposition \ref{S3subgroup} we need two more lemmas.

\begin{lemma}\label{equation}
Let $A\in Gl(n,\z\p)$. Then, there exists a vector 
$v=(v_1,...,v_n)\in (\z\p)^n$ such that $v_i$ is a  
square of a $p$-adic unit for all $i$ 
and such that $Av$ is a vector whose components are given by elements of 
$\z\lp$.
\end{lemma}

\begin{proof}
Let $B:=A^{-1}$. We have to solve the following problem:
Find a vector $w\in \z\lp^n$ such that $Bw\neq 0$ has as components
squares of $p$-adic units. 
The question whether a $p$-adic unit is a 
square, can be decided by reducing to $\z/p$ for $p$ odd or to $\z/8$ for 
$p=2$.
In both cases the reduction $\quer B$ of $B$ is an 
invertible matrix and induces 
therefore an epimorphism on $(\z/p)^n=:V$. In particular, if $\quer v \in V$ 
is a vector with components given by squares mod $p$ such that 
all entries are units in $\z/p$, 
there exists a vector
$w\in (\Z)^n$ such that $\quer B w=v$. Hence, $Bw$ is a vector whose components
are squares of nontrivial $p$-adic units. For $p=2$, 
the same argument works, we only have to replace $\z/p$ by $\z/8$.
\end{proof}

{\it Proof of Proposition \ref{S3subgroup}.}
Let $M$ and $N$ denote the loop spaces and $G$ the Lie group 
constructed in  Proposition \ref{globalsubgroup}
Since $BM\p\simeq BG\p$ and $BM_0\simeq BG_0$
we have a pull back diagram
$$
\diagram
BM\lp \rrto \dto & & BG\p \dto \\
BG_0 \rto & BG^\flat_p \rto^A & (BG\p)_0
\enddiagram
$$
Here $BG^\flat_p$ is the formal p-adic completion of the rationalization
$BG_0$ in the sense of Sullivan, and $(BG\p)_0$ the localization at 0 of  
$BG\p$. The map $A$ is an equivalence between the homotopy equivalent spaces 
$BG^\flat_p$ and $(BG\p)_0$, and induces a continuous map in homotopy.
The homotopy groups $\pi_*((BG\p)_0$ carry a natural topology, since 
$\pi_*(BG\p) \cong \pi_*(BG)\otimes \z\p$ (details may be found in 
\cite{wilkerson}).
$(BG\p)_0\cong K(\q\p{}^r,4)$ is a rational Eilenberg-MacLane space. 
Since self maps of rational Eilenberg-MacLane spaces 
are determined by the induced maps in homotopy,
and since $A$ induces a continuous map
in homotopy, we can think of $A$ as a matrix in $Gl(n,\q\p)$ inducing a 
continuous self equivalence of $(\q\p)^n$. 
Such matrices can be written as a product $BR$ 
where $B\in GL(n,\z\p)$ and $R \in Gl(n;\q)$. For example, this follows 
from the fact that the adic genus of products of $S^1$'s is rigid.
Since $R$ can be realized as a 
self equivalence of $BG_0$, replacing $A$ by $B$ does not change 
the homotopy type of the pull back. Hence we may assume that 
$A\in Gl(n,\z\p)$. 

We have an analogous pull back diagram as above for the classifying space 
$B\widetilde M$ of the universal cover $\widetilde M$ of $M$ with the same 
gluing map $A$, namely
$$
\diagram
B\widetilde M\lp \rrto \dto & & BS^3{}\p\times B\widetilde H\p \dto \\
BS^3_0\times B\widetilde H_0 \rto & 
(BS^3\times B\widetilde H)^\flat_p \rto^A & ((BS^3\times B\widetilde H)\p)_0
\enddiagram
$$ 
Here we used the fact that $G\cong S^3\times H$ with $H = \widetilde H/K$, 
and $\widetilde H$ a product of $S^3$'s. Since every equivalence
$BS^3{}\p \lra BS^3{}\p$ induces in $\pi_4(BS^3{}\p)$ multiplication 
by a non trivial square unit of 
$\z\p$, Lemma \ref{equation}
shows that there exists 
a map $BS^3 \lra BS^3{}\p \times B\widetilde H\p$ 
such that
the composition  
$$
BS^3\lp \lra BS^3{}\p \lra BS^3{}\p \times B\widetilde H\p \lra 
((BS^3\times B\widetilde H)\p)_0
\larrow{A^{-1}} (BS^3\times B\widetilde H)^\flat_p
$$
lifts to a map $BS^3\lp \lra BS^3_0\times B\widetilde H_0$. Moreover,
localized at 0, composition with the projection on the first factor 
is an equivalence.
This establishes a map $BS^3\lp \lra B\widetilde M\lp$ such that 
the completion of the composite $BS^3\lp \lra B\widetilde M\lp \lra BM\lp$ 
is induced by the 
monomorphism
$S^3{}\p \lra S^3{}\p \times \widetilde H\p \lra 
S^3{}\p\times \widetilde H/K\p = G\p$ of 
p-compact groups.
This shows that the homotopy fiber of 
$BS^3\lp \lra BM\lp$ is simple and $\z\lp$-finite
as is the homotopy fiber of the 
composite $f : BS^3\lp \lra BM\lp \larrow{g} BN\lp$.
Since $H^4(BN_0;\q) \cong H^4(BM_0;\q)$, there exists a left inverse
$s: BN_0 \lra BM_0$ for $g_0$. Projection onto the first factor
gives a left inverse of 
$BS^3_0 \lra BM_0 \simeq BS^3_0\times B\widetilde H_0$. This shows that, 
localized at 0,
the map $f : BS^3\lp \lra BN\lp$ has a left inverse and finishes the proof of 
the proposition.
\qed
\end{section}

\begin{section}{Stable Reducibility of abelian quotients}

Let $L$ be a finite loop space. Let $T$ be an abelian compact Lie group, 
and let $B\varphi : BT \rightarrow BL$ be a monomorphism. We will show in 
this section that the $\z$-finite fiber $L/T$ is stably reducible.
Unfortunately, here we use a stronger assumption than we are able to produce
in our case. In general, the fibration sequence
$S^1 \lra L \lra L/S^1 \lra BS^1$ cannot be extended one further step 
to the right. This is only possible after completion. 

We begin 
with the case of a Lie group $L$ to develop our intuition.

\begin{lemma}\label{Liecase}
Let $\varphi : T \rightarrow L$ be a monomorphism between compact Lie groups, 
where $T$ is abelian, then $L/T$ is stably reducible.
\end{lemma}
\begin{proof}
Let $\call$ and $\calt$ be the Lie algebras of $L$ and $T$ respectively. 
It is easy to see that the tangent bundle of $L/T$ is given by 
$L \times_T (\call/\calt)$, where the action of $T$ on the Lie algebras 
is given by the adjoint representation. Since $T$ is abelian, the bundle 
$L \times_T \calt$ is trivial. On adding it to the tangent bundle of $L/T$, 
we get the bundle $L \times_T \call$. Notice that the adjoint action of $T$ 
on $\call$ extends to $L$ and hence $L \times_T \call$ is trivial. This shows 
that $L/T$ has a stably trivial tangent bundle, which is equivalent to 
 being stably reducible.
\end{proof}

One would like to extend this argument to the case of $L$ being a finite 
loop space. The theory that best preserves the analogy with compact Lie 
groups has been developed by John Klein in \cite{Klein}. For a 
topological group $G$, Klein defines the dualizing $G$-spectrum, 
$D_G = \map(EG_+,S[G])^G$, where $S[G]$ is the spectrum 
$\Sigma^{\infty} G_+$ with the left $G$ action. The residual 
right $G$ action on $S[G]$ induces the $G$ action on $D_G$. The 
spectrum $D_G$ is the appropriate notion of the Adjoint representation. 
This is justified by the fact that for a compact Lie group $G$, there is 
an equivalence $D_G = S^{Ad}$, where $S^{Ad}$ denotes the one-point 
compactification of the adjoint representation of $G$. 

The magic of the spectrum $D_G$ lies in the following two theorems of 
Klein \cite{Klein}

\begin{theorem}\label{PD}
Assume that $BG$ is a finitely dominated space. Then the following are 
equivalent:
\nz
(i) $BG$ is a Poincar\'{e} duality space of formal dimension $n$.
\nz
(ii) $D_G$ has the (unequivariant) homotopy type of a sphere spectrum of 
dimension $-n$.

Moreover, in the above two cases, the Thom spectrum $EG_+\smash_G D_G$ 
is the Thom spectrum of the Spivak normal bundle of $BG$.
\end{theorem}

\begin{theorem}\label{ext}
Assume $ 1\rightarrow H \rightarrow G \rightarrow Q \rightarrow 1$ is 
an extension of topological groups, then if $BH$ is a finitely 
dominated Poincar\'{e} duality space, then there is a weak equivalence 
of spectra
\[ D_G \cong D_H \smash D_Q \]
Moreover, one may replace $D_H$ by an $H$-equivalent spectrum so as to 
make the above equivalence $H$-equivariant. 
\end{theorem}

We will apply the above theorems to the extension of topological groups
\[ 1 \rightarrow \Omega(L/T) \rightarrow T \rightarrow L \rightarrow 1 \]
In order to make sense of this extension, we must work in the 
model-category of simplicial groups and replace the above groups by 
equivalent topological groups (c.f \cite{Go-Ja}).

Let us record a simple lemma about simplicial groups that will be useful 
in the sequel.

\begin{lemma}\label{fib}
Let $sH \rightarrow sG$ be an acyclic fibration of simplicial groups. 
Then on taking realizations one gets an extension of topological groups
\[ 1 \rightarrow K \rightarrow H \rightarrow G \rightarrow 1 \]
with a contractible kernel $K$. Moreover, there is an $H$-equivalence of 
spectra $D_H \cong D_G$. 
\end{lemma}
\begin{proof}
The realization of an acyclic fibration of simplicial groups is an 
acyclic Serre fibration. Hence we get an extension of topological groups 
with a contractible kernel $K$. We may consider the space $EH/K$ as a 
model for $EG$. Then the required equivalence is induced by 
the (2-sided) $H$-equivalence $S[H] \rightarrow S[G]$ and is given by 
\[ 
D_H = \map(EH_+,S[H])^H \rightarrow \map(EH_+,S[G])^H = 
\map(EG_+,S[G])^G = D_G 
\]
\end{proof}

We are now ready to prove the main theorem of this section
\begin{theorem}\label{red}
Let $L$ be a finite loop space of formal dimension $n$, and let 
$B\varphi : BT \rightarrow BL$ be a monomorphism, where $T$ is an 
abelian compact Lie group of rank $r$, then the $\z$-finite space 
$L/T$ is stably reducible.
\end{theorem}
\begin{proof}
For a connected space $X$, let $s\Omega X$ denote the simplicial Kan 
loop group of $X$ (c.f \cite{Go-Ja}). The map $B\varphi$ induces a 
simplicial homomorphism
\[ s\varphi : s\Omega BT \rightarrow s\Omega BL \]
The properties of model categories allow us to factor $s\varphi$ 
through an acyclic cofibration followed by a fibration. Consequently, 
we may assume that $s\varphi$ is a fibration. Unfortunately, 
$s\Omega BT$ is a highly non-abelian model for $T$. This is the price 
one has to pay for obtaining a fibration. Fortunately however, 
the simplicial group $s\Omega BT$ is related to an abelian 
simplicial model for $T$ via adjointness. We have an acyclic fibration  
\[ s\Omega BT \rightarrow sT^a \]
where $sT^a$ is a simplicial abelian model for $T$. 

Now taking the realization of $s\varphi$, we obtain an extension
\[ 
1 \rightarrow K \rightarrow T^f \rightarrow L^f \rightarrow 1\ , 
\]
where $T^f$ and $L^f$ denote the free models of $T$ and $L$, we obtain 
by realizing the simplicial Kan loop groups $s\Omega BT$ and $s\Omega BL$ 
respectively. The group $K$ is clearly equivalent to $\Omega(L/T)$. 
Using Theorem \ref{ext}, one obtains a $K$-equivariant equivalence 
\[ 
D_{T^f} \cong D_K \smash D_{L^f} 
\]
Hence $K$ acts trivially on $D_{L^f}$. In fact, $D_{L^f}$ is a 
sphere spectrum of dimension $n$ with a trivial $K$ action. One 
sees this as follows: by Theorem \ref{PD}, $D_K$ is a sphere of 
dimension $r-n$. Moreover, by Lemma \ref{fib} 
$D_{T^f} \cong D_{T^a} \cong D_{T}$ is also a sphere of dimension $r$, 
thus $D_{L^f}$ is a sphere of dimension $n$ with a trivial $K$ action. 
Hence $D_K \cong \Sigma^{-n}D_{T^f}$ as a $K$-spectrum. It now follows 
from Theorem \ref{PD} that the Thom spectrum of the Spivak Normal bundle 
of $L/T$ is given by $\Sigma^{-n}EK_+\smash_K D_{T^f}$. Thus to show 
that $L/T$ is stably reducible, it is sufficient to show that 
$D_{T^f}$ is equivalent to a sphere spectrum with a trivial 
$T^f$ action. This follows by applying Lemma \ref{fib} to notice 
that $D_{T^f}$ is equivalent to $D_{T^a}$ which clearly has a trivial 
$T^f$ action. 
\end{proof}  

\end{section}
\bigskip

\end{document}